\input amstex
\documentstyle{amsppt}
\magnification=\magstep1

\define\eps{\varepsilon}
\define\BH{\Bbb B(\Cal H)}

\overfullrule=0pt
\TagsOnRight

\topmatter

\title {A note on fixed points of \\
completely positive maps}\endtitle

\rightheadtext{Fixed points of completely positive maps}

\author{Gert K. Pedersen} \endauthor

\date{July  2003}\enddate

\address {Department of Mathematics, University of Copenhagen, 
Universitetsparken 5, DK-2100 Copenhagen \O, Denmark.}
\endaddress

\email {gkped\@math.ku.dk}\endemail


\endtopmatter

\document

\footnote""{\copyright 2003 by the author. This paper may be
reproduced, in its entirety, for non-commercial purposes.}

Given a Hilbert space $\Cal H$ let $(x_t)_{t\in T}$ be a 
strongly (equivalently weakly) measurable field of operators 
in $\BH$ over a locally compact Hausdorff space $T$. Assume 
moreover that $\int_t x_t^*x_t\, d\mu(t)=\bold 1$ and 
$\int_T x_tx^*_t\, d\mu(t)=e\le \bold 1$ for some Borel 
measure $\mu$ on $T$, where the integrals are computed, say, 
by considering the corresponding sesquilinear forms, e.g.
$(e\xi|\eta) = \int_T (x_t^*\xi|x_t^*\eta)\,d\mu(t)$
for $\xi,\eta$ in $\Cal H$.

This means in particular that the linear map 
$$
\Phi\colon\BH\to\BH \quad\text{given by}\quad
\Phi(a)=\int_T x_t^*ax_t\, d\mu(t)
$$ 
is a completely positive, unital map. Usually completely 
positive maps are presented by a sequence, i\.e\. by 
taking $T=\Bbb N$ and $\mu$ the counting measure, but 
there are cases where the continuous model may be preferred. 
Note that the condition $\int x_tx_t^*\,d\mu(t)\le \bold 1$ 
puts a restraint on the type of completely positive maps to 
be described. It is most obviously satisfied if $x_t=x_t^*$
for all $t$ in $T$.

If $a$ is an element in $\BH$ commuting with the field 
$(x_t)_{t\in T}$, then evidently $\Phi(a)=a$. The surprising 
fact is that in many cases this is the only way fixed points 
can arise.

If $M$ is a semi-finite von Neumann subalgebra of $\BH$ we
say that an element $a$ in $M_+$ is {\it finite}, if 
$\tau_i(a)<\infty$ for a separating family $\{\tau_i\}$ of 
normal, semi-finite traces on $M$.

\proclaim{Theorem} With $\Phi$ and $M$ as above, assume 
that $\Phi(M)\subset M$ and that $a$ is a finite element
in $M_+$ with $\Phi(a)\ge a$. Then $\Phi(a)=a$ and as a 
consequence $ax_t=x_ta$ for $\mu-$almost all $t$ in $T$.
\endproclaim

\demo{Proof} Let $\tau$ be a normal, semi-finite trace on $M$
with $\tau(a)<\infty$. Then
$$
\aligned
&\tau(\Phi(a))=\int_T \tau(x_t^*ax_t)\, d\mu(t)
=\int_T \tau(a^{1/2}x_tx_t^*a^{1/2})\,d\mu(t) \\
=\;&\tau\left( a^{1/2}\left(\int_T x_tx_t^*\,d\mu(t)\right)a^{1/2}\right)
=\tau(a^{1/2}ea^{1/2})\le\tau(a)\,.
\endaligned
\tag{$*$}
$$
It follows that $\Phi(a)$ is a finite element of $M_+$. If now
$\Phi(a)\ge a$ then $\Phi(a)-a$ is a finite element in $M_+$ and 
by ($*$)
$$
0\le\tau(\Phi(a)-a)=\tau(\Phi(a))-\tau (a)\le 0.
$$
Thus $\Phi(a)-a$ is annihilated by a separating family of 
semi-finite, normal traces on $M$, whence $\Phi(a)-a=0$, as 
desired. This proves the first statement in the Theorem.

For each real $\eps$ define $f_\eps (t)=t^2(1-\eps t)^{-1}$ on 
the open interval $]-|\eps|^{-1}, |\eps|^{-1}[$. It follows 
from the results of Bendat and Sherman, \cite{\bf 2} or 
\cite{\bf 6}, that $f_\eps$ is an operator convex function. As 
shown in \cite{\bf 6}, see also \cite{\bf 7}, this implies 
that $f_\eps$ satisfies the {\it Jensen operator inequality}
$$
f_\eps\left(\sum_{k=1}^n x_k^*ax_k\right) 
\le \sum_{k=1}^n x_k^*f_{\eps}(a)x_k
$$
for every finite set of operators $(x_k)$ with $\sum x_k^*x_k \le \bold 1$.
Since $f_\eps$ is operator continuous on any closed, bounded 
subinterval of $]-|\eps|^{-1}, |\eps|^{-1}[$ it follows by standard
approximation arguments that we then also have 
$$
f_\eps\left(\int_T  x_t^*ax_t\,d\mu(t)\right) 
\le \int_T x_t^*f_{\eps}(a)x_t\,d\mu(t)
$$
whenever $(x_t)_{t\in T}$ is a measurable operator field with 
$\int_T x_t^*x_t\,d\mu(t)\le \bold 1$.

In our case this means that 
$$
f_{\eps}(a)=f_\eps(\Phi(a))\le \Phi(f_\eps(a)).
$$
Since $f_\eps(a)$ (for $|\eps|<\Vert a\Vert^{-1}$) is a finite 
element in $M_+$ (dominated by $\lambda a$, where $\lambda 
=\Vert a\Vert (1-|\eps|\Vert a\Vert)^{-1}$), it follows from the 
first part of the proof that actually $f_{\eps}(a)=\Phi(f_\eps(a))$.
Expanding $f_\eps$ in a norm convergent series we therefore have
$$
a^2+\eps a^3+\eps^2 a^4 + \ldots 
= \Phi(a^2)+\eps\Phi(a^3)+\eps^2\Phi(a^4) + \ldots
$$
for every $\eps$ in a neighbourhood of zero, from which we conclude 
that $\Phi(a^n)=a^n$ for every $n$.

It follows that $\Phi(f(a))=f(a)$ for every polynomial $f$, hence 
by continuity and approximation for every bounded measurable 
function $f$. In particular $\Phi(p)=p$ for every spectral 
projection $p$ of $a$. Arguing as in \cite{{\bf 4}, Lemma 3.3} this 
means that 
$$
0=(\bold 1-p)p(\bold 1-p)=\int_T (\bold 1-p)x_t^*px_t(\bold 1-p)\,d\mu(t),
$$
whence $px_t(\bold 1-p)=0$ for almost all $t$ in $T$. However,
this computation is also valid for the spectral projection 
$\bold 1-p$, so $(\bold 1-p)x_tp=0$ as well. We conclude that $x_tp=px_t$
almost everywhere, and since this holds for every spectral 
projection, also $x_ta=ax_t$ almost everywhere, as desired.
\hfill$\square$
\enddemo
 
\bigskip

\proclaim{Corollary} If $\Phi(a)=a$ and $a^2$ is finite in $M$,
then $x_ta=ax_t$ for $\mu-$almost all $t$ in $T$.
\endproclaim

\demo{Proof} Since $\Phi(a)^2\le\Phi(a^2)$ for every completely 
positive, unital map, the condition $\Phi(a)=a$ implies that 
$a^2\le\Phi(a^2)$. If now $a^2$ is finite we must have 
$\Phi(a^2)=a^2$ by the Theorem. Consequently $x_tp=px_t$ for
$\mu-$almost all $t$ in $T$ and every spectral projection $p$
for $a^2$. Since $a$ and $a^2$ have the same spectral family
we conclude that $x_ta=ax_t$ almost everywhere.
\hfill$\square$
\enddemo

\bigskip

\example{Remarks} With $T=\Bbb N$, $x_n=x_n^*$ and $M=\BH$ 
(so that $a$ is a trace class operator) the result in the 
Theorem is known to phycisists as the L\"u{}ders Theorem, 
cf\. \cite{\bf 10} and \cite{\bf 3}. It is shown in \cite{{\bf 1},
Theorem 4.2} that the result may fail in general, i\.e\. without
extra conditions on either $a$ or $(x_n)$. Indeed, if the 
von Neumann subalgebra $A$ of $\BH$ generated by the family 
$(x_n)$ is not injective then there is an operator $a$ such 
that $\Phi(a)=a$, but $a\notin A'$, cf\. \cite{{\bf 1}, Theorem 3.6}.

That the condition af finiteness on $a$, although sufficient, is 
not necessary is shown in \cite{{\bf 3}, Proposition 1} and also 
in \cite{{\bf 1}, Theorem 3.2}. Indeed, if $a$ is an element in 
$\BH_+$ with a pure point spectrum that can be totally ordered in 
decreasing order, and if furthermore $x_t=x_t^*$ for all $t$ in $T$, 
then $\Phi(a)=a$ (or just $\Phi(a)\ge a$) implies that $ax_t=x_ta$ 
for $\mu-$almost all $t$ in $T$. The simple argument is reproduced 
below.

By assumption there is an orthogonal family $(p_n)$ of projections 
with sum {\bf 1} such that $a=\sum \lambda_n p_n$ (strongly 
convergent sum), where $\lambda_1 > \lambda_2 >\ldots \ge 0$. 
If $\Phi(a)\ge a$ then for each unit vector $\xi$ in the eigenspace
$p_1(\Cal H)$ we have
$$
\lambda_1 = (a\xi|\xi)\le (\Phi(a)\xi|\xi)
=\int_T (ax_t\xi|x_t\xi)\,d\mu(t)
\le \lambda_1 \int_T (x_t\xi|x_t\xi)\,d\mu(t)=\lambda_1,
$$
since $\Vert a\Vert =\lambda_1$. Consequently 
$\lambda_1\bold 1 - a \ge 0$ and
$$
\int_T ((\lambda_1\bold 1-a)x_t\xi|x_t\xi)\,d\mu(t)=\lambda_1-\lambda_1=0,
$$
whence $ax_t\xi=\lambda_1x_t\xi$ for $\mu-$almost all $t$ in $T$.
Thus $x_tp_1(\Cal H)\subset p_1(\Cal H)$, which implies that 
$p_1x_t=x_tp_1$ since $x_t=x_t^*$. Using that $\Phi(p_1)=p_1$ 
we pass to the operator $a_1=a-\lambda_1p_1$ in $\BH_+$, which 
again satisfies the condition $\Phi(a_1)\ge a_1$, and proceed by 
induction to show that $p_nx_t=x_tp_n$ for $\mu-$almost all 
$t$ and every $n$, whence $ax_t=x_ta$, as desired.

This result, which contains the original L\"u{}ders Theorem, shows
that $\Phi(a)=a$ implies that $a\in \{x_n\}'$ for every positive,
compact operator. Making the generalizing from $\BH$ to a 
semi-finite factor $M$ of type II with faithful, semi-finite trace 
$\tau$, one could therefore hope that the Theorem would be valid 
(when $x_t=x_t^*$) not only for elements $a$ in $M_+$ with 
$\tau(a)<\infty$ or $\tau(a^2)<\infty$, but for all positive 
elements in the norm closed ideal $M_\tau$ generated by the finite 
elements (or just the finite projections) in $M$.
\endexample

\enddocument